\documentclass[oneside]{amsart}
\newcommand{\formatswitch}{preprint}

\usepackage{amssymb}
\usepackage{ifthen}
\usepackage{verbatim}
\usepackage{psfig}
\usepackage[all]{xy}

\newcommand{\tref}[1]{(\ref{#1})}

\ifthenelse{\equal{\formatswitch}{big}}{
\setlength{\textwidth}{432pt}
\setlength{\oddsidemargin}{18.8775pt}

\setcounter{tocdepth}{1}
}{
\ifthenelse{\equal{\formatswitch}{paper}}{

\setcounter{tocdepth}{1}
}{
\setcounter{tocdepth}{2}
}
}
\DeclareMathAlphabet\EuScript{U}{eus}{m}{n}
\DeclareMathAlphabet\EuScriptb{U}{eus}{b}{n}

\newcommand{\claimenum}{\renewcommand{\theenumi}{\alph{enumi}}
 \renewcommand{\labelenumi}{\textit{(\theenumi)}}
 \renewcommand{\theenumii}{\roman{enumii}}
 \renewcommand{\labelenumii}{\textit{(\theenumii)}}
 \begin{enumerate}}
\newcommand{\claimenumend}{\end{enumerate}}

\newcommand{\romanenum}{\renewcommand{\theenumi}{\roman{enumi}}
 \renewcommand{\labelenumi}{\textit{(\theenumi)}}
 \renewcommand{\theenumii}{\alph{enumii}}
 \renewcommand{\labelenumii}{\textit{(\theenumii)}}
 \begin{enumerate}}
\newcommand{\romanenumend}{\end{enumerate}}

\newtheorem{dummy}{realdumb}[section]
\newtheorem{thm}{Theorem}

\newtheorem{lemma}[dummy]{Lemma}
\newtheorem{prop}[dummy]{Proposition}
{\theoremstyle{definition} \newtheorem{defn}[dummy]{Definition}}

{\theoremstyle{definition} }
\newtheorem{cor}{Corollary}[dummy]

\renewcommand{\text}{\mathrm}

\newcommand{\strutdepth}{\dp\strutbox}
\newcommand{\marginalnote}[1]
   {\strut\vadjust{\kern-\strutdepth\domarginalnote{#1}}}
\newcommand{\domarginalnote}[1]{\vtop to \strutdepth{
  \baselineskip\strutdepth
   \vss\llap{ #1\ \ }\null}}  
\newcounter{showlabelflag}
\newcounter{makelabelflag}
\newcommand{\showlabels}{\setcounter{showlabelflag}{1}}

\newcommand{\makelabels}{\setcounter{makelabelflag}{1}}
\newcommand{\hidelabels}{\setcounter{showlabelflag}{2}}
\newcommand{\mylabel}[1]{
  \ifthenelse{\value{makelabelflag}=1}
    {\label{#1}}{}
  \ifthenelse{\value{showlabelflag}=1}
    {\marginpar{#1}}{}\relax}

\newcommand{\Z}{{\mathbf Z}}
\newcommand{\N}{{\mathbf N}}
\newcommand{\sub}{\subseteq}

\ifthenelse{\equal{\formatswitch}{draft}}{
\showlabels
}{\relax}

\newtheorem{question}{Question}
\newtheorem{conj}[question]{Conjecture}

\newcommand{\LZ}{\ensuremath {}^\infty(\wr\,\Z)}
\newcommand{\RZ}{\ensuremath (\Z\,\wr)^\infty}

\newcommand{\symm}{\text{Sym}}

\newcommand{\mymargin}[1]{
  \ifthenelse{\value{showlabelflag}=1}
    {\marginpar{#1}}{}\relax}

\newcounter{enumo}

\newcommand{\RRsh}{\kern -1 pt \Rsh}

\newcounter{keepitemnum}

\newcounter{keepitemnumm}

\begin{document}

\bibliographystyle{amsplain}
\title[Elementary Amenable Subgroups]{Elementary Amenable Subgroups \\
of R.~Thompson's group \(F\)}
\thanks{AMS
Classification (2000): primary 20E07, secondary 20F19, 20F22, 20E22.}
\author{MATTHEW G. BRIN}
\date{June 29, 2004}
\CompileMatrices


\makelabels
\hidelabels
\maketitle


\section{Introduction}

This paper contributes to the study of the subgroups of Thompson's
group \(F\) by constructing a sequence of subgroups of increasing
complexity.  The group \(F\) is an interesting finitely presented
group with a pleasant, faithful representation in the group
\(PL_o(I)\) of orientation preserving, piecewise linear, self
homeomorphisms of the unit interval.  Thus our subgroups also lie in
\(PL_o(I)\).  The subgroups we construct are all elementary amenable
of a certain bounded degree (see Section \ref{TheClassesSec} for
definitions) while \(F\) is not elementary amenable.  Thus the
complexities of the subgroups that we construct are in some sense
strictly less than that of \(F\) itself and we raise the question of
whether there are subgroups of \(PL_o(I)\) that are in the gap
between the groups we construct and the group \(F\).  We also raise
(and give specific meaning to) the question of whether \(F\) is the
``only'' subgroup of \(PL_o(I)\) that is not elementary amenable.
See Conjecture \ref{AltConj}.

The group \(F\) originated with Thompson in his studies in algebraic
logic \cite{thomp:notes}, was independently discovered as the key
behind interesting examples in shape and homotopy theory
\cite{dydak}, \cite{dydak2}, \cite{freyd+heller}, can be called the
``structure group'' of the associative law \cite{dehornoy:assoc},
and arises in a fundamental way in the theory of diagram groups
\cite{guba+sapir2}.  It is also the source of difficult questions,
the most prominent of which is whether \(F\) is amenable
\cite{gersten:problems}.  See \cite{CFP} and \cite[Section
4]{brown:finiteprop} for an introduction to Thompson's groups.

\subsection{Background}

The subgroup structure of \(F\) is only partly understood.  Some
information about subgroups of \(F\) comes from the faithful
representation of \(F\) in \(PL_o(I)\).  For example, \(PL_o(I)\) is
torsion free and no subgroup of \(PL_o(I)\) can be free on two
generators (\cite{brin+squier} or \cite{CFP}).  Further, every
non-abelian subgroup of \(PL_o(I)\) contains a copy of \(\Z\wr\Z\),
the wreath product of \(\Z\) with itself \cite[Theorem
21]{guba+sapir:diag2}.  Other information comes from the fact that
\(F\) is (one of the prime examples of) a diagram group
\cite{guba+sapir2}.  For example, every nilpotent subgroup of \(F\)
is abelian \cite[Cor. 15]{guba+sapir:diag2} .  (One can also show
that every nilpotent subgroup of \(PL_o(I)\) is abelian.)

Other results that will be reflected in the themes of this paper are
that \(F\) contains many copies of itself as well as many wreath
products \cite{brin:ubiq}, and in particular contains the iterated
wreath product \((\cdots((\Z\wr\Z)\wr\Z)\wr\cdots)\)
\cite[Cor. 20]{guba+sapir:diag2}.  Not relevant to this paper are
various results about the metric structure of subgroups of \(F\).
See \cite[Section 6]{guba+sapir:diag2}, \cite{burillo},
\cite{MR2001k:20087} and \cite{MR2016187}.

\subsection{The classes \protect\(EG_\alpha\protect\),
\protect\(EG\protect\), \protect\(AG\protect\) and
\protect\(NF\protect\)}\mylabel{TheClassesSec}

Thompson's group and its subgroups are scattered among several
classes of groups from \cite{MR19:1067c}: the class \(EG\) of
elementary amenable groups, the class \(AG\) of amenable groups, and
the class \(NF\) of groups with no subgroup isomorphic to the free
group on two generators.  The class \(NF\) needs no further
definitions, the definitions relevant to the class \(EG\) will be
given below, and the definitions relevant to the class \(AG\) are to
be found in \cite {MR19:1067c}, \cite{MR87e:04007} or
\cite{MR94g:04005} as well as full discussions on the relations
between the classes.  It is known (\cite{MR19:1067c},
\cite{MR94g:04005}) that \(EG\sub AG\sub NF\), that both
containments are proper (\cite{MR82b:43002}, \cite{MR84m:43001} and
\cite{grig2}), and that there are finitely presented groups in each
difference (\cite{MR99b:20055} and \cite{MR2004f:20061}).

Let \(EG_0\) be the class of all groups that are either finite or
abelian.  The class \(EG\) is defined in \cite{MR19:1067c} as the
smallest class of groups that contains \(EG_0\) and is closed under
the four operations of (I) taking subgroups, (II) taking quotients,
(III) taking extensions and (IV) taking direct limits.  It is shown 
in \cite{chou} operations (I) and (II) are not necessary.  We need
some details of this.

A hierarchy \(EG_\alpha\) of classes of groups indexed by the
ordinals is defined inductively in \cite{chou} as follows.  With
\(EG_0\) as defined above, let \(\alpha>0\) be an ordinal so that
\(EG_\beta\) is defined for all \(\beta<\alpha\).  If \(\alpha\) is
not a limit ordinal, then \(EG_{\alpha}\) is the class of all groups
\(G\) that are extensions \[1\rightarrow N\rightarrow G\rightarrow
Q\rightarrow 1\] with \(N\) and \(Q\) in some \(EG_\beta\) with
\(\beta<\alpha\), or are the direct limits of groups from the various
\(EG_\beta\) with \(\beta<\alpha\).  If \(\alpha\) is a limit
ordinal, then \(EG_\alpha\) is the union of all the \(EG_\beta\)
with \(\beta<\alpha\).  The following combines Propositions 2.1 and
2.2 of \cite{chou}.

\begin{prop}[Chou]\mylabel{chou} The following are true of the items
defined above.  {\begin{enumerate}
\renewcommand{\theenumi}{\alph{enumi}} \item Each \(EG_\alpha\) is
closed under operations (I) and (II). \item \(EG\) is the union of
the \(EG_\alpha\).  \item \(EG\) is the smallest class of groups
that contains \(EG_0\) and is closed under operations (III) and
(IV).  \end{enumerate}} \end{prop}

It is remarked in \cite{chou} that since the class of isomorphism
types of finitely generated groups has the cardinality of the
continuum, there must be some \(EG_{\delta}\) containing all
finitely generated groups in \(EG\), and since every group is the
direct limit of its finitely generated subgroups, we must have
\(EG=EG_{\delta+1}\).

If \(G\) is a group in \(EG\), then the {\itshape elementary class}
of \(G\) (often referred to in this paper as the {\itshape class} of
\(G\)) is the smallest ordinal \(\alpha\) so that \(G\) is in
\(EG_\alpha\).

We use standard notation for certain ordinals.  The ordinal
\(\omega\) is the smallest infinite ordinal, setting
\(1\omega=\omega\) gives the definition of \((n+1)\omega\) for
\(n\ge1\) as the limit of \(n\omega, n\omega+1, n\omega+2, \ldots\),
and \(\omega^2\) as the limit of \(1\omega, 2\omega, 3\omega,
\ldots\).

\subsection{Results and questions}

The main result of this paper is the following.

\begin{thm}\mylabel{MainThm} For each non-limit ordinal \(\alpha\le
\omega^2+1\), there is a subgroup of \(F\) that is of elementary
class \(\alpha\).  \end{thm}

We obtain Theorem \ref{MainThm} by finding techniques for building
groups of higher class from groups of a given class or set of
classes.  We observe that these techniques work well as long as we
are working primarily with finitely generated groups and that our
technique runs out of finitely generated groups at class
\(\omega^2+1\).  The question arises as to whether it is possible to
get past \(\omega^2+1\) in either \(F\) or \(PL_o(I)\).

It is remarked above that \(PL_o(I)\) is in the class \(NF\) and it
is known \cite{CFP} that \(F\) is not in \(EG\).  It is very easy
\cite{brin:ubiq} for isomorphic copies of \(F\) to appear in
subgroups of \(PL_o(I)\).  While not stated as such, our
constructions are designed to increase elementary class while
avoiding the inclusion of copies of \(F\).

We now formally bundle the various observations that we have made
above, the theorems of this paper, and the results of
\cite{brin:ubiq} and \cite{guba+sapir:diag2} into various questions
and conjectures.

The observation of \cite{guba+sapir:diag2} that every non-abelian
subgroup of \(PL_o(I)\) contains a copy of \(\Z\wr\Z\) has the form
of an ``alternative'' theorem: every subgroup of \(PL_o(I)\) is
either abelian or contains \(\Z\wr\Z\) as a subgroup.  The group
\((\cdots((\Z\wr\Z)\wr\Z)\wr\cdots)\), which we will denote by
\(\RZ\), is not solvable, while all finite iterations of
the wreath product of copies of \(\Z\) are solvable.  Sapir has
raised the question of whether every non-solvable subgroup of \(F\)
contains a copy of \(\RZ\).  The question can also be
asked of \(PL_o(I)\).  Our second result concerns another
non-solvable wreath product \[\LZ=
(\cdots\wr(\Z\wr(\Z\wr\Z))\cdots)\] where the wreath product is that
of permutation groups (see Section \ref{ConstrSec} for definitions)
and not the standard (restricted) wreath product of
\cite{neumann:wreath}.

\begin{thm}\mylabel{NoCrossEmbds} Both \(\LZ\) and
\(\RZ\) are isomorphic to subgroups of \(F\).  Further,
there is no subgroup of \(\RZ\) that is isomorphic to 
\(\LZ\), and there is no subgroup  of
\(\LZ\) that is isomorphic to  \(\RZ\).  \end{thm}

From Theorem \ref{NoCrossEmbds}, the answer to Sapir's question is
``No.''  However, the following can be asked (which we put in
``alternative'' form).

\begin{question}\mylabel{SolvAlt} Is it true of every subgroup of
\(PL_o(I)\) that it is either solvable or it contains a subgroup
that is isomorphic to one of \(\LZ\) or \(\RZ\)?  \end{question}

The next question relates to our failure to find elementary amenable
subgroups of class greater than \(\omega^2+1\).

\begin{question}\mylabel{HowHighTheMoon} Is \(\omega^2+1\) the
largest elementary class of an elementary amenable subgroup of
\(PL_o(I)\)?  \end{question}

The following has been privately suggested by the author.

\begin{conj}\mylabel{AltConj} Every subgroup of \(PL_o(I)\) is
either elementary amenable or it contains a subgroup that is
isomorphic to \(F\).  \end{conj}

If we raise Question \ref{HowHighTheMoon} to the status of
conjecture and combine it with  Conjecture \ref{AltConj}, we get the
following.

\begin{conj}\mylabel{BigAltConj} Every subgroup of \(PL_o(I)\) is
either elementary amenable of some class \(\alpha\) with \(\alpha\le
\omega^2+1\), or it contains a subgroup that is isomorphic to \(F\).
\end{conj}

We include a remark that has little to do with this paper.

The group \(F\) is contained in \(NF\smallsetminus EG\), but it is
not known whether it is contained in \(AG\).  The group \(F\) has a
presentation that is somewhat smaller than that of the example in
\cite{MR99b:20055} of a group in \(AG\smallsetminus EG\) and that is
vastly smaller than that of the example in \cite{MR2004f:20061} of a
group in \(NF\smallsetminus AG\).  Thus it is still of interest to
locate \(F\) accurately in the string of classes \(EG\sub AG\sub
NF\).  It has also been put forth by Grigorchuk that answering the
question of the amenability of \(F\) will be an important step in
understanding the groups in the class \(NF\).

\subsection{The Thompson group \protect\(V\protect\)}

Thompson's group \(V\) (see \cite{CFP}) is more flexible in that it
acts on the totally discontinuous Cantor set.  In particular it
contains copies of the free group on two generators.  It will be
clear to those familiar with \(V\), that there is much more room to
apply the constructions of this paper in the setting of \(V\).  It
would be of interest to know more about the subgroups of \(V\).
That this is not a trivial task is indicated by the fact that \(V\)
contains a copy of the infinite symmetric group and thus a copy of
every finite group.  It would also be of interest to know if there
is a criterion along the lines of \cite{brin:ubiq} to recognize
\(V\) among the homeomorphisms of the Cantor set.

\subsection{Section contents}

In Section \ref{ConstrSec} we develop construction techniques that
will apply to the setting of \(PL_o(I)\).  Our main tool will be
the permutation wreath product and we set up methods for recognizing
the product.  Section \ref{ElemFactsSec} discusses methods for
recognizing the elementary class of groups with emphasis on
recognizing the elementary class of wreath products.  Section
\ref{MainThmSec} proves Theorem \ref{MainThm} while Section
\ref{ZInfinitySec} proves Theorem \ref{NoCrossEmbds}.

In the narrative below, lemmas without proofs are to be taken as
exercises. 

\subsection{Thanks} The author would like to thank Fernando Guzm\'an
for numerous helpful conversations.

\section{Construction techniques}\mylabel{ConstrSec}

\subsection{Permutation groups}

A {\itshape permutation group} is a pair \((G,Z)\) in which \(Z\) is
a set and \(G\) is a subgroup of \(\symm(Z)\), the symmetric group
on \(Z\) or the group of all bijections from \(Z\) to itself.  We
have no need of a more general notion.

The group \(G\) will act on \(Z\) on the right and \(xG=\{xg\mid
g\in G\}\) is the {\itshape orbit} of \(x\) under \(G\).  The
{\itshape support} of \(G\) written \(\text{Supp}(G)\) is the set
\(\{x\mid \{x\}\ne xG\}\).  There is also an action of \(G\) on the
subsets of \(Z\) in the usual way: if \(X\) is a subset of \(Z\),
then \(Xg\) denotes the image of \(X\) under the action of \(g\in
G\) and is the set \(\{xg\mid x\in X\}\).  Note that this makes
\(XG\) a collection of sets.

\subsection{The permutation wreath
product}\mylabel{PermWreathDefSec}

We need to recognize the wreath product of two permutation groups.
If \((G,A)\) and \((H,B)\) are permutation groups in which the
action of \(H\) on \(B\) is transitive, then the following
definition of the permutation group \((G\wr H, A\times B)\) is
equivalent to the definition found in \cite[Section 1.6]{drobinson}.
We define \(G\wr H\) as the subgroup of \( \symm(A\times B)\)
generated by \(G'\) and \(H^*\) where \(G'\) and \(H^*\) are images
of two homomorphisms \(G\rightarrow G'\) and \(H\rightarrow H^*\).
The first homomorphism is defined by choosing a fixed \(b_1\in B\)
and sending \(g\in G\) to \(g'\in G'\) where for and \((a,b)\in
A\times B\), \[(a,b)g' =\begin{cases} (ag,b), &b=b_1, \\(a,b), &b\ne
b_1, \end{cases}\] and the second is defined by sending \(h\in H\)
to \(h^*\in H^*\), where \((a,b)h^*=(a,bh)\).  The definition does
not depend on the choice of \(b_1\).

In the event that the action of \(H\) on \(B\) is not transitive,
then the permutation wreath product as defined in \cite{drobinson}
is not generated by \(H^*\) and \(G'\), but instead by \(H^*\) and
an isomorphic copy of \(G\) for each orbit in \(B\) under the action
of \(H\).  In our setting, we do not have to worry about this.  Note
that we do not require that the action of \(G\) on \(A\) be
transitive.

The reader can glean the properties of the permutation wreath
product from \cite{drobinson} including the associativity of the
product which is reflected in Corollary \ref{DoubleWreath} below.

\subsection{Pre-wreath structures}

The wreath product involves two factors.  To make it easier to
iterate the product, we give conditions under which it is easy to be
a factor.

\begin{defn}\mylabel{PreWreathDef} A {\itshape pre-wreath structure}
is a quadruple \((Z, Y, H, X)\) in which \(H\) is a non-trivial
group and the rest are sets and the following are satisfied.
\begin{enumerate} \item \(H\sub \symm(Z)\).  \item
\(\text{Supp}(H)\sub Y\sub Z\).  \item \(\emptyset\ne X\sub Y\).
\item For all \(h\in H\), we have that \(Xh\cap X\ne \emptyset\)
implies \(h|_X=1|_X\).  \item For all \(1\ne h\in H\), there is a
\(j\in H\) so that \(Xjh\ne Xj\).  \end{enumerate} \end{defn}

We have no need of trivial groups here and the non-triviality
assumption on \(H\) will be convenient.  The need for having both
\(Y\) and \(Z\) will not be apparent until Proposition
\ref{PreWreathToWreath} below.

The collection \(\{Xh\mid h\in H\}\) will be important and we will
denote it by \(XH\).  We will use \(\overline{XH}\) to denote the
union \(\displaystyle{\bigcup_{h\in H}Xh}\) of the elements in
\(XH\).

\begin{lemma}\mylabel{PreFacts}  Let \((Z, Y, H, X)\) be a
pre-wreath structure.  Then the following hold.
{\begin{enumerate} 
\renewcommand{\theenumi}{\alph{enumi}}
\item \(XH\) is a collection of pairwise disjoint non-empty subsets
of \(Y\).  
\item \(X\in XH\).  
\item For each \(h\in H\) and \(A\in XH\) we have \(Ah\in XH\).  
\item Each non-identity element \(h\in H\) has some \(A\in XH\) so
that \(Ah\ne A\).  
\item For every \(A\) and \(B\) in \(XH\), there is an \(h\in H\)
with \(Ah=B\).  
\item For each \(h, k\in H\) and \(A\in XH\), if \(Ah=Ak\), then
\(h|_A=k|_A\).  
\end{enumerate}}
\end{lemma}

From (c), there is a homomorphism \(H\rightarrow \symm(XH)\).  Item
(d) implies that this homomorphism is an injection.  In particular
\((H, XH)\) is a permutation group if we identify \(H\) with its
image under the injection \(H\rightarrow \symm(XH)\).

Other observations are that (e) implies that the action of \(H\) on
\(XH\) is transitive, and (f) implies that for any \(A\) and \(B\)
in \(XH\), there is only one function from \(A\) to \(B\) that is a
restriction of an element of \(H\).

In a pre-wreath structure \((Z, Y, H, X)\), we call \(XH\) the
{\itshape carrier} of the structure.

\begin{lemma}\mylabel{SubPreWreath} If \((Z, Y, H, X)\) is a
pre-wreath structure, and if \(X'\) is a non-empty subset of \(X\),
then \((Z, Y, H, X')\) is a pre-wreath structure.  \end{lemma}

\subsection{Pre-wreath structures to wreath products}

We need a standard fact.  

\begin{lemma}\mylabel{PushProblemsLeft} If \(G\) and \(H\) are
subgroups of a group \(M\), then any element of \(\langle G,
H\rangle\) can be written as an element of \(H\) followed by a
product of elements of \(G\) each of which is conjugated by an
element of \(H\).  \end{lemma}

\begin{proof} We have
\[\begin{split} &h_ng_{n-1} \cdots h_3g_2h_2g_{1}h_1 = \\
&\quad
(h_n\cdots h_3h_2h_1) \\
&\qquad
(h_1^{-1}h_2^{-1}h_3^{-1}\cdots h_{n-1}^{-1} g_{n-1} h_{n-1}\cdots
h_3h_2h_1)
\\
&\quad\qquad 
\cdots \\
&\qquad\qquad 
(h_1^{-1}h_2^{-1}h_3^{-1}g_3 h_3h_2 h_1)  \\
&\qquad\qquad\quad
(h_1^{-1}h_2^{-1}g_2h_2h_1) \\
&\qquad\qquad\qquad
(h_1^{-1}g_1h_1).
\end{split} \]
\end{proof}

\begin{prop}\mylabel{PreWreathToWreath} Let \((Z, Y, H, X)\) and
\((Z, X, G, W)\) be pre-wreath structures.  Then \[\left(Z, Y,
\langle G, H\rangle, \, W\right)\] is a pre-wreath structure and its
carrier \(W\langle G, H\rangle\) is equal to \[WGH=\{Wgh\mid g\in
G,\,\,h\in H\}.\]

Further, the permutation group \((\langle G, H\rangle, WGH)\) is
isomorphic to (similar to) the permutation group \((G\wr H,
(WG)\times (XH))\) as defined in \cite{drobinson}.  \end{prop}

The picture below might help understand the hypotheses of
Proposition \ref{PreWreathToWreath}.  We only indicate \(Y\) and
some of its subsets since \(Y\) contains all of the ``action.''

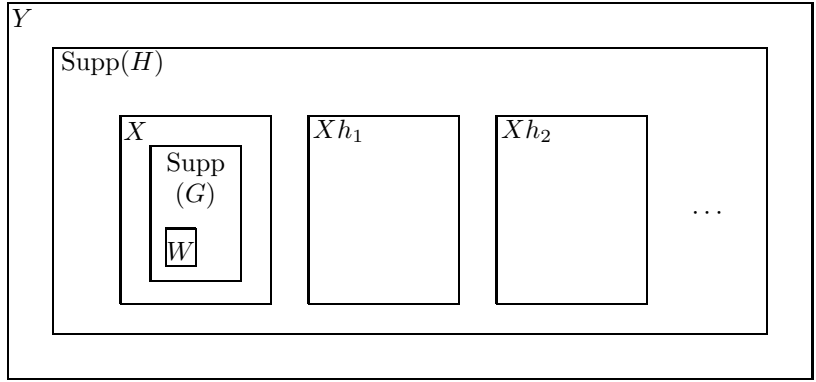
\[
\xy
(5,5); (112,5)**@{-}; (112,55)**@{-}; (5,55)**@{-};
(5,5)**@{-};
(7,53)*{Y};
(11,11); (106,11)**@{-}; (106,49)**@{-}; (11,49)**@{-};
(11,11)**@{-};
(19,47)*{\text{Supp}(H)};
(20,15); (40,15)**@{-}; (40,40)**@{-}; (20,40)**@{-}; (20,15)**@{-};
(22,38)*{X}; 
(24,18); (36,18)**@{-}; (36,36)**@{-}; (24,36)**@{-}; (24,18)**@{-};
(30,33.5)*{\text{Supp}}; (30,29.5)*{(G)};
(26,20); (30,20)**@{-}; (30,25)**@{-}; (26,25)**@{-}; (26,20)**@{-};
(28,22)*{W};
(45,15); (65,15)**@{-}; (65,40)**@{-}; (45,40)**@{-}; (45,15)**@{-};
(49,38)*{Xh_1}; 
(70,15); (90,15)**@{-}; (90,40)**@{-}; (70,40)**@{-}; (70,15)**@{-};
(74,38)*{Xh_2}; 
(98,27)*{\cdots};
\endxy
\]

In the following proof, (1)--(5) refer to items in Definition
\ref{PreWreathDef} and (a)--(f) refer to items in Lemma
\ref{PreFacts}.

\begin{proof} Denote \(\langle G, H\rangle\) by \(M\).  The group
\(M\) is contained in \( \symm(Z)\) with support in \(Y\).  We have
\(\emptyset\ne W\sub X\sub Y\).  Thus (1), (2) and (3) are obtained
with no problem.

For the rest of the proof, we need a better understanding of the
group \(M\).  Let \(w\in M\).  By Lemma \ref{PushProblemsLeft},
\(w\) is an element of \(H\) followed by a product of elements of
\(G\) conjugated by elements of \(H\).

Note that the conjugate \(G^h=h^{-1}Gh\) has support in \(Xh\).
Because of (4) applied to \((Z, Y, H, X)\), we know that the
conjugate \(G^h\) is completely determined by the set \(Xh\).  Since
the elements \(Xh\) of \(XH\) are pairwise disjoint, the group
\[K=\langle G^h \mid h\in H\rangle\] is the direct sum of the
elements of \(\{G^h\mid h\in H\}\).  By (f) applied to \((Z, Y, H,
X)\), the action of \(H\) on \(K\) is to permute the direct factors
and we get that \(M\) is the semidirect product \(K\rtimes H\).

Consider \(Ww\) with \(w\in M\).  We know \(w=hv\) with \(h\in H\)
and \(v\in K\), so \(Ww=Whv\sub Xhv\).  All components of \(v\) in
the product structure of \(K\) fix \(Xh\) except the component
coming from \(G^h\) which has support in \(Xh\).  Thus \(Ww\sub
Xh\).  

If \(Ww\cap Ww'\ne\emptyset\) with \(w=hv\) and \(w'=h'v'\), then
\(Xh\cap Xh'\ne\emptyset\) and \(h|_W=h'|_W\) since \(h|_X=h'|_X\)
and \(W\sub X\).  Now the actions of \(v\) and \(v'\) on \(Xh=Xh'\)
are determined by their components from \(G^h=G^{h'}\).  Let
\(h^{-1}gh\) and \(h^{-1}g'h\) be these components from \(v\) and
\(v'\), respectively.  It follows that \(Wg\cap Wg'\ne\emptyset\),
so \(g|_W=g'|_W\).  Since \(W\sub X\) and \(h|_X=h'|_X\), we get
that \((hg)|_W=(h'g')|_W\).

Let \(1\ne w=hv\) be from \(M\).  If \(h\ne 1\), then there is a
\(j\in M\) so that \(Xjh\ne Xj\).  Now \(Wj\sub Xj\) has \(Wjw\sub
Xjhv=Xjh\) which is disjoint from \(Xj\).  If \(h=1\), then
\(v\ne1\) and some component \(l^{-1}gl\) of \(v\) is not trivial.
There is a \(k\) from \(G\) with \(Wkg\ne Wk\).  Now
\(Wklw=Wklv=Wkl(l^{-1}gl) =Wkgl\ne Wkl\) since \(l\) is a bijection.
This verifies that \((Z, Y, M, W)\) is a pre-wreath structure.

That \(WGH\sub W\langle G, H\rangle=WM\) is clear.  If \(w\in M\) is
of the form \(w=hv\) with \(h\in H\) and \(v\in K\), then let
\(h^{-1}gh\) be the component of \(v\) from \(G^h\).  Now
\(Ww=Whv=Wh(h^{-1}gh) = Wgh\) and we have \(WM\sub WGH\).

To establish the claimed similarity in the last statement, we note
that sending \((Wg, Xh)\in (WG)\times (XH)\) to \(Wgh\sub Xh\) is
clearly a surjection from \((WG)\times (XH)\) onto \(WGH\).  It is
one-to-one since there is only one bijection from \(X\) to \(Xh\)
that is a restriction of an element of \(H\).  Use \(t\) to denote
the bijection \((Wg, Xh)\mapsto Wgh\).

We note that \(M\) is defined as \(M = \langle G, H\rangle\) in \(
\symm(Z)\).  Since (d) holds for the pre-wreath structure \((Z, Y,
M, W)\), we know that the image of \(M\) in \( \symm(WGH)\) is
isomorphic to \(M\).  Since \(M\) is generated by \(G\) and \(H\),
its image in \( \symm(WGH)\) is generated by the images of \(G\) and
\(H\).  We will now use \(M\), \(G\) and \(H\) to denote their
images in \( \symm(WGH)\).

The group \(G\wr H\) is the subgroup of \( \symm((WG)\times (XH))\)
generated by \(G'\) and \(H^*\) where \(g\in G\) is sent to \(g'\in
G'\) and \(h\in H\) is sent to \(h^*\in H^*\), and where for and
\((Wg_1,Xh_1)\in (WG)\times (XH)\), \[(Wg_1,Xh_1)g' =\begin{cases}
(Wg_1g,Xh_1), &Xh_1=X, \\(Wg_1,Xh_1), &Xh_1\ne X \end{cases}\] and
\((Wg_1,Xh_1)h^*=(Wg_1,Xh_1h)\).

Now \[\begin{split}(Wg_1, Xh_1)tgt^{-1} = (Wg_1h_1)gt^{-1} =
&\begin{cases} (Wg_1g)t^{-1}, &Xh_1=X, \\ (Wg_1h_1)t^{-1}, & Xh_1\ne
X, \end{cases}\\ = &\begin{cases}(Wg_1g,X), &Xh_1=X, \\ (Wg_1,
Xh_1), &Xh_1\ne X. \end{cases}\end{split}\] Also
\[(Wg_1,Xh_1)tht^{-1}=(Wg_1h_1)ht^{-1}=(Wg_1h_1h)t^{-1}=(Wg_1,
Xh_1h)\] and we see that \(t\) conjugates the actions of \(G\) and
\(H\) to those of \(G'\) and \(H^*\).  This establishes the
similarity and completes the proof.  \end{proof}

\begin{cor} \mylabel{DoubleWreath} Let \[\begin{split} &(Z, Y_1,
H_1, Y_0), \\ &(Z, Y_2, H_2, Y_1), \\ &(Z, Y_3, H_3, Y_2)
\end{split}\] be pre-wreath structures and let \(A=Y_0H_1\times
Y_1H_2\times Y_2H_3\).  Then \[(Z, Y_3, \langle H_1, H_2,
H_3\rangle, Y_0)\] is a pre-wreath structure with carrier
\[T=Y_0\langle H_1, H_2\rangle H_3 = Y_0\langle H_1, H_2, H_3\rangle
= Y_0H_1\langle H_2, H_3\rangle = Y_0H_1H_2H_3\] and the permutation
group \((\langle H_1, H_2, H_3\rangle, T)\) is isomorphic to
(similar to) \[((H_1\wr H_2)\wr H_3, A) \simeq (H_1\wr (H_2\wr H_3),
A).\] \end{cor}

In Proposition \ref{PreWreathToWreath}, the status of \(\langle G,
H\rangle\) as a wreath product needs much less than was assumed.  We
leave it to the reader to extract a proof of the following from the
proof of Proposition \ref{PreWreathToWreath}.

\begin{lemma}\mylabel{PreWreathToWreathII} Let \((Z, Y, H, X)\) be a
pre-wreath structure, and let \((G,Z)\) be a permutation group with
support in \(X\).  Then restricting the action of \(\langle G,
H\rangle\) to \(\overline{XH}\) gives a permutation group \((\langle
G, H\rangle, \overline{XH})\) that is isomorphic to (similar to) the
permutation group \((G\wr H, X\times (XH))\).  \end{lemma}

\subsection{Homomorphisms and normal subgroups} There is much known
about normal subgroups in wreath products \cite{neumann:wreath}.
We need very little and give below what we need.  We start with a
useful homomorphism.

\begin{lemma}\mylabel{TopGpHom} In the setting of Lemma
\ref{PreWreathToWreathII}, there is a homomorphism from \(\langle G,
H\rangle\) to \(H\) that is the identity on \(H\) and whose kernel
is the group \(K\) generated by the conjugates of \(G\) in \(\{G^h
\mid h\in H\}\).  \end{lemma}

\begin{proof} By Lemma \ref{PushProblemsLeft}, any element \(w\) of
\(\langle G, H\rangle\) can be written as an element \(h\) of \(H\)
multiplied by conjugates of elements of \(G\).  The conjugates of
elements of \(G\) have supports in the elements of \(XH\), so \(w\)
and \(h\) agree on \(Z\smallsetminus \overline{HX}\) and affect the
same permutation on the sets of \(XH\).  But Lemma \ref{PreFacts}(f)
says that an element of \(H\) is determined by its permutation on
\(XH\) and its action on \(Z\smallsetminus \overline{HX}\).  Thus
sending \(w\) to \(h\) is well defined.  It is clearly a
homomorphism and its kernel is \(K\).  \end{proof}

The next lemma gives some control over the kernels and quotients.

\begin{lemma}\mylabel{KernelControl} Assume the setting of Lemma
\ref{PreWreathToWreathII} and let \(N\) be a normal subgroup of
\(\langle G, H\rangle\).  
{\begin{enumerate}
\renewcommand{\theenumi}{\Roman{enumi}}
\item If \(N\) contains an element whose permutation on \(XH\) is
non-trivial, then there is a subgroup of \(N\) that surjects onto
\(G\). 
\item If every element of \(N\) has trivial permutation on \(XH\),
then \(\langle G, H\rangle/N\) surjects onto \(H\).
\end{enumerate}}  \end{lemma}

\begin{proof} Under the assumptions of (I), the transitivity of the
action of \(H\) on \(XH\) gives an \(f\) in \(N\) so that \(Xf\ne
X\), implying that \(Xf\cap X=\emptyset\).  If \(g\) is any element
of \(G\), then \(f^{-1}g^{-1}fg\) is an element of \(N\), its
support is in \(X\cup Xf\), its restriction to \(X\) equals that of
\(g\), and its restriction to \(Xf\) is that of \(f^{-1}g^{-1}f\).
This is a non-homomorphic injection of \(G\) into \(G\times G^f\)
whose image we will call \(\overline{G}\).  The subgroup of
\(G\times G^f\) generated by \(\overline{G}\) is a subgroup of \(N\)
whose support is in \(X\cup Xf\) and whose restriction to \(X\) is
\(G\).

Under the assumptions of (II), it follows from Lemma
\ref{PushProblemsLeft} and the injection of \(H\) into the symmetric
group on \(XH\), that every element of \(N\) is a product of
conjugates of \(G\).  Thus \(N\) is contained in the kernel \(K\) of
Lemma \ref{TopGpHom} and the result follows from that lemma.
\end{proof}

\subsection{An infinite construction}

Let \((Z, Y_i, H_i, Y_{i-1})\) be pre-wreath structures for positive
integers \(i\).  Let \[H'_i = \langle H_1, H_2, H_3, \cdots,
H_i\rangle.\] We have \(H'_i \sub H'_{i+1}\) and we can let \(H\) be
the union of all the \(H'_i\).  We also have \(Y_i\sub Y_{i+1}\) and
we can let \(Y\) be the union of all the \(Y_i\).  Let
\(T'_i=Y_0H_1H_2H_3\cdots H_i\).  Inductively from Proposition
\ref{PreWreathToWreath} and its corollary, each \((Z, Y_i, H'_i,
Y_0)\) is a pre-wreath structure, \(T'_i\) is its carrier and
\(T'_i=Y_0H'_i\).  Clearly, \(T'_i\sub T'_{i+1}\) and we let \(T\)
be the union of all the \(T'_i\).

\begin{lemma}\mylabel{FirstIndConstr}The quadruple \((Z, Y,
H, Y_0)\) is a pre-wreath structure with carrier \(T\).  \end{lemma}

\begin{proof} We check (1)--(5) in Definition \ref{PreWreathDef} to
show that \(( Z, Y, H,Y_0)\) is a pre-wreath structure.  Items (1)
and (3) are immediate.  For (2), we note that an element of \(H\) is
in one of the \(H'_i\) and has its support in \(Y_i\sub Y\).  Items
(4) and (5) are also immediate from the fact that each element of
\(H\) is in one of the \(H'_i\), and \((Z, Y_i, H'_i, Y_0)\) is a
pre-wreath structure by an inductive extension of Corollary
\ref{DoubleWreath}.

The carrier of \((Z, Y, H, Y_0)\) is \(Y_0H=\{Y_0h\mid h\in H\}\).
But a given \(h\in H\) is in some \(H'_i\) and \(Y_0h\) is in
\(Y_0H'_i\).  But \(Y_0H'_i=T'_i\).  This gives \(Y_0H\sub
T\).  However, each element of \(T'_i\) is clearly in \(Y_0H\), so
\(T\sub Y_0H\).  \end{proof}

\subsection{Isomorphisms}

If \(j:Z\rightarrow Z'\) is a bijection of sets, then \(j\) induces
an isomorphism \(j_*:\symm(Z)\rightarrow \symm(Z')\) by
``conjugation'' in that \(z'(\sigma j_*) = z'(j^{-1}\sigma j)\).
The quotes can be removed from the word ``conjugation'' if we put
\(Z\) and \(Z'\) in the disjoint union of \(Z\) and \(Z'\).

An isomorphism \(f:(Z, Y, H, X)\rightarrow (Z', Y', H', X')\) of
pre-wreath structures is a bijection \(f:Z\rightarrow Z'\) so that
\(Y'=Yf\), \(X'=Xf\) and \(H'=Hf_*\).  We could add the isomorphism
\(f_*|H:H\rightarrow H'\) to the data of an isomorphism of a
pre-wreath isomorphism, but it is not necessary.

\begin{lemma} {\upshape (I)} If \(f:(Z, Y, H, X)\rightarrow
(Z', Y', H', X')\) is an isomorphism of pre-wreath structures, then
\((XH)f=(X'H')\).  That is, \(f\) takes carrier onto carrier.

{\upshape (II)} If \((Z, Y, H, X)\) is a pre-wreath structure, and
\(f:Z\rightarrow Z'\) is a bijection, then \((Z', Yf, Hf_*, Xf)\) is
a pre-wreath structure and \(f:(Z, Y, H, X)\rightarrow (Z, Yf, Hf_*,
Xf)\) is an isomorphism. \end{lemma}

\subsection{An inductive situation}

Let \((Z, Y_1, H_1, Y_0)\) be a pre-wreath structure, and let
\(f:Z\rightarrow Z\) be a bijection so that \(Y_0f=Y_1\).  Define
\[\begin{alignedat}{2} Y_i &=Y_0f^i=Y_1f^{i-1}, &&i\in\Z \\ H_i &=
H_1 (f_*)^{i-1}, &\qquad&i\in\Z.  \end{alignedat}\] This gives that
each \((Z, Y_i, H_i, Y_{i-1})\) is a pre-wreath structure for each
\(i\in\Z\) with carrier \(T_i=Y_{i-1}H_i\).  We have \(T_i=T_{i-1}f
= T_1f^{i-1}\) for \(i\in\Z\).

Recall that \(\overline{T}_i\) is the union of all the sets in
\(T_i\).  We also have \(\overline{T}_i = \overline{T}_{i-1}f =
\overline{T}_1f^{i-1}\) for \(i\in\Z\).

We now have a doubly ended sequence of pre-wreath structures and can
create various singly ended sequences just by picking different
starting places.  Thus we have a sequence of singly ended sequences
to which Lemma \ref{FirstIndConstr} applies.  We establish notation
for the pre-wreath structures that result from applying that lemma.

For \(j\in\Z\), let \((Z, \widetilde{Y}_j, R_j, Y_{j-1})\) be the
pre-wreath structure that results when Lemma \ref{FirstIndConstr} is
applied to the sequence of pre-wreath structures \((Z, Y_i, H_i,
Y_{i-1})\) for \(i\ge j\).  The notation \(R_j\) is used since the
group is the result of an iteration of wreath products \(R_j =
H_j\wr H_{j+1}\wr\cdots\) that goes off to the right.  The
associativity of the permutation wreath product makes parentheses
unnecessary.

It is our intention to make a pre-wreath structure based on the
union of the groups \(R_j\).  However we need to make one slight
adjustment in the situation so that a carrier of this union is easy
to construct.

We know that \(\text{Supp}(H_1)\sub Y_1\).  From now on we are going
to assume  
\mymargin{IndAssumpt} \begin{equation}\tag{6} \label{IndAssumpt}
Y_0\sub \text{Supp}(H_1)\sub Y_1 \quad \text{and}\quad
\text{Supp}(H_1)\ne Y_1.\end{equation} Since \(Y_0\sub
\text{Supp}(H_1)\) and \(Y_0f=Y_1\), we have a non-empty subset
\(W_0\) of \(Y_0\) so that \(W_0f\cap \text{Supp}(H_1)=\emptyset\).
Thus \(W_0f\sub (Y_1\smallsetminus \text{Supp}(H_1))\).  If we
define \(W_i=W_0f^i\) for \(i\in\Z\), then it follows from
\(Y_i=Y_{i-1}f\) and \(\text{Supp}(H_i)=\text{Supp}(H_{i-1})f\) for
all \(i\in\Z\) that \(W_i\sub(Y_i\smallsetminus \text{Supp}(H_i))\).
The following picture might help.

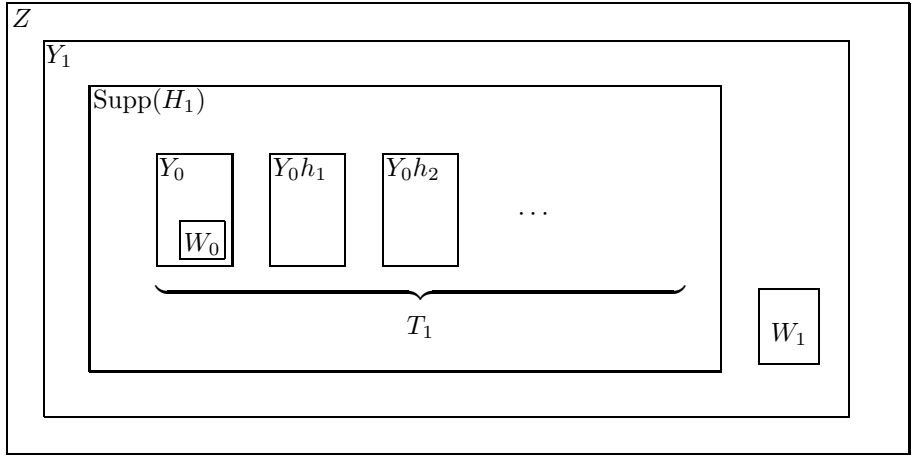
\[
\xy
(0,0); (120,0)**@{-}; (120,60)**@{-}; (0,60)**@{-}; (0,0)**@{-};
(2,58)*{Z};
(5,5); (112,5)**@{-}; (112,55)**@{-}; (5,55)**@{-};
(5,5)**@{-};
(7,53)*{Y_1};
(11,11); (95,11)**@{-}; (95,49)**@{-}; (11,49)**@{-};
(11,11)**@{-};
(19,47)*{\text{Supp}(H_1)};
(20,25); (30,25)**@{-}; (30,40)**@{-}; (20,40)**@{-}; (20,25)**@{-};
(22,38)*{Y_0}; 
(23,26); (29,26)**@{-}; (29,31)**@{-}; (23,31)**@{-}; (23,26)**@{-};
(26,28)*{W_0};
(35,25); (45,25)**@{-}; (45,40)**@{-}; (35,40)**@{-}; (35,25)**@{-};
(39,38)*{Y_0h_1}; 
(50,25); (60,25)**@{-}; (60,40)**@{-}; (50,40)**@{-}; (50,25)**@{-};
(54,38)*{Y_0h_2}; 
(70,32)*{\cdots};
(55,22)*{\underbrace{
  \xy
  (0,0)*{\hbox to1pt{}}; 
  (70,0)*{\hbox to1pt{}};
  \endxy}};
(55,17)*{T_1};
(100,12); (108,12)**@{-}; (108,22)**@{-}; (100,22)**@{-}; (100,12)**@{-};
(104,16)*{W_1}; 
\endxy
\]

A trivial conclusion from the fact that \(H_1\) is not trivial and
\(Y_0\) is not empty is that \(\text{Supp}(H_1)\) is not empty and
contains all of \(Y_0H_1\).  In turn, the non-empty
\(\text{Supp}(H_1)\) is contained in \(Y_1=Y_0f\) which is contained
in \(\text{Supp}(H_2)=(\text{Supp}(H_1))f\).  Eventually,
\(\text{Supp}(H_i)\sub \text{Supp}(H_j)\) for all \(i<j\).

\begin{lemma}\mylabel{SecondIndConstr} Let \((Z, Y_1, H_1, Y_0)\) be
a pre-wreath structure, and let \(f:Z\rightarrow Z\) be a bijection
so that \(Y_0f=Y_1\).  Define \(Y_i\), \(H_i\), \(T_i\),
\(\widetilde{Y}_j\), and \(R_j\) as above.  Assume that
\tref{IndAssumpt} holds, let \(\emptyset\ne W_0\sub Y_0\) so that
\(W_0f\cap \text{Supp}(H_1)=\emptyset\), and let \(W_i=W_0f^i\) for
\(i\in \Z\).  Then for each \(j\in \Z\), the following are true.
{\begin{enumerate}
\renewcommand{\theenumi}{\alph{enumi}}
\item \((Z, \widetilde{Y}_j, R_j, Y_{j-1})\) and \((Z,
\widetilde{Y}_j, R_j, W_{j-1})\) are pre-wreath structures.
\item \(f:(Z, \widetilde{Y}_j, R_j, Y_{j-1})\rightarrow (Z,
\widetilde{Y}_{j+1}, R_{j+1}, Y_{j})\) and \(f:(Z,
\widetilde{Y}_j, R_j, W_{j-1})\rightarrow (Z,
\widetilde{Y}_{j+1}, R_{j+1}, W_{j})\) are isomorphisms.
\item \(R_j\sub R_{j+1}\).
\item The various carriers are related by \(W_{j-1}R_j \sub
Y_{j-1}R_j\sub Y_jR_{j+1}\), and \(\overline{Y_{j-1}R_j}\,\,\cap
\,\, \overline{W_jR_{j+1}} = \emptyset\).
\item For all \(k\ne j\), the underlying sets
\(\overline{W_{j-1}R_j}\) and \(\overline{W_{k-1}R_k}\) are
disjoint.
\item For \(j<k\) the sets in \(W_{k-1}R_{k}\) are disjoint from
\(\text{Supp}(H_j)\).  
\end{enumerate}}
\end{lemma}

\begin{proof} Item (a) follows from Lemmas \ref{FirstIndConstr} and
\ref{SubPreWreath}.  The verifications for (b) and (c) are
elementary.  For (d), the containments follow from \(W_0\sub Y_0\sub
Y_1\) and the properties of the action of \(f\) and its powers.  If
we show that \(\overline{Y_0R_1}\,\,\cap\,\, \overline{W_1R_2} =
\emptyset\), then the disjointness claim in (d) will then follow by
applying powers of \(f\).

By Lemma \ref{FirstIndConstr}, \[Y_0R_1 = \bigcup_{i\ge1}
Y_0H_1H_2\cdots H_i = \bigcup_{i\ge2}T_1H_2H_3\cdots H_i.\]  A
similar formula for \(Y_1R_2\) and 
Definition \ref{PreWreathDef}(4) gives 
\[W_1R_2=\bigcup_{i\ge2} W_1H_2H_3\cdots H_i = \bigcup_{i\ge2}
W_0fH_2H_3\cdots H_i.\]
It suffices to show
\(\overline{W_0fH_2H_3\cdots H_i}\,\,\cap\,\, \overline{T_1
H_2H_3\cdots H_i} = \emptyset\) for each \(i\ge2\).

We already know that \(W_0f\cap \overline{T}_1=\emptyset\).  Both
\(W_0fH_2H_3\cdots H_i\) and \(\overline{T}_1 H_2H_3\cdots H_i\) are
collections of subsets of \(Y_i\).  If we have disjoint subsets
\(A\) and \(B\) of \(Y_i\), then we are done if we show that
\(AH_{i+1}\) and \(BH_{i+1}\) are disjoint.  But if they are not,
then \(Ah\cap Bh'\ne\emptyset\) for some \(h\) and \(h'\) in
\(H_{i+1}\), so that \(Y_ih\cap Y_ih'\ne\emptyset\).  This implies
that \(h|_{Y_i}=h'|_{Y_i}\) which is impossible if the disjoint
\(A\) and \(B\) in \(Y_i\) have intersecting images under \(h\) and
\(h'\).

Item (e) follows from the truth of (d) for all \(j\in\Z\).

Lastly, \(\text{Supp}(H_{k-1})\) is disjoint from \(W_{k-1}\) and
both are in \(Y_{k-1}\).  By (4) of Definition \ref{PreWreathDef},
every element of \(R_k\) that takes points of \(Y_{k-1}\) into
\(Y_{k-1}\) is the identity on \(Y_{k-1}\).  Thus
\(\text{Supp}(H_{k-1})\) is disjoint from all sets in
\(W_{k-1}R_k\).  Now every \(\text{Supp}(H_j)\) is contained in
\(\text{Supp}(H_{k-1})\) for \(j<k\).  \end{proof}

\begin{lemma}\mylabel{ThirdIndConstr}  Assume the hypotheses and
notation of Lemma \ref{SecondIndConstr}.  Let \(L_i\) be the group
\(\langle H_{i-1}, H_{i-2}, \ldots\rangle\), and let \(M\) be the
union of all the \(R_i\).  Then \(M\) is also the union of all the
\(L_i\) and is also generated by the union of all the \(H_i\).
Further, for each \(i\in\Z\), the following are true.
{\begin{enumerate}
\renewcommand{\theenumi}{\alph{enumi}}
\item \(f:L_i\rightarrow L_{i+1}\) is an isomorphism.
\item  \(M=L_i\wr R_i\).  
\item For each \(k>0\), \(L_{i+k}=L_i\wr \langle H_i, H_{i+1},
\ldots, H_{i+k-1}\rangle\).
\item Each non-trivial normal subgroup of \(L_i\) and each
non-trivial normal subgroup of \(M\) contains a subgroup that
surjects onto \(L_0\).
\end{enumerate}}
\end{lemma}

\begin{proof} That statements about what constitutes \(M\) are
immediate.  Item (a) is immediate.  Item (b) follows from the fact
that all of the \(H_j\) are represented in \(L_i\cup R_i\), the fact
that \((Z, \widetilde{Y}_i, R_i, Y_{i-1})\) is a pre-wreath
structure, the fact that the support of \(L_i\) is in \(Y_{i-1}\),
and from Lemma \ref{PreWreathToWreathII}.  Item (c) is similar.  To
argue (d), we note that a non-trivial normal subgroup \(N\) of
\(L_i\) or \(M\) contains an element that moves a set in some
\(W_{j-1}R_j\).  Now by Lemma \ref{KernelControl} applied to
\(M=L_j\wr R_j\) or \(L_i=L_j\wr\langle H_{j+1}, \ldots,
H_{i-1}\rangle\), a subgroup of \(N\) surjects onto \(L_j\) which is
isomorphic to \(L_0\).  \end{proof}

Note that we do not claim a pre-wreath structure for \(M\) or for
any of the \(L_i\).

\begin{prop}\mylabel{FourthIndConstr}  Assume the hypotheses and
notation of Lemma \ref{SecondIndConstr} and let 
\[V
=\langle M, f\rangle
=\langle L_1, f\rangle
=\langle R_1, f\rangle
=\langle H_1, f\rangle.\]
Then \(V\) is the ascending HNN extension of \(R_1\) given by the
injection of \(R_1\) in itself under conjugation by \(f\) and it is
also the semidirect product of \(M\) with the infinite cyclic group
generated by \(f\). Further, \((Z,
{Y}_1\cup\text{Supp}(f), V, W_0)\) is a pre-wreath
structure whose carrier \(W_0V\) is the union of the \(W_0R_1f^n\),
\(n\in\Z\).  Lastly, any non-trivial normal subgroup of \(V\)
contains a subgroup that surjects onto \(L_0\).  \end{prop} 

\begin{proof} The claimed equalities are immediate and the claimed
structures as HNN extensions and semidirect products follow from the
standard characterizations of those structures.

We look at the five items in Definition \ref{PreWreathDef} of a
pre-wreath structure.  Item (1) is immediate.  Item (2) follows from
the fact that the support of \(H_1\) is in \(Y_1\).  Item (3) is
immediate.

To work on the remaining items and the claim about \(W_0V\), we
consider an arbitrary element \(v\) of \(V\).  From Lemma
\ref{PushProblemsLeft}, \(v=f^n h_1h_2h_3\cdots h_k\) where each
\(h_i\) is in a conjugate of \(H_1\) by a power of \(f\).  Thus for
each \(i\) there is a \(j_i\) so that \(h_i\in H_{j_i}\).  Write
\(v\) as \(v=f^nw\) with \(w=h_1h_2h_3\cdots h_k\).  Writing
\(M=L_{n+1}\wr R_{n+1}\), let \(\phi_n:M\rightarrow R_{n+1}\) be as
given by Lemma \ref{TopGpHom} and let \(\overline{w}\) be the image
of \(w\) under \(\phi_n\).  This is just \(w\) with all \(h_i\)
removed that have \(j_i\le n\).

Now \(W_0v = W_0f^n w\) and \(w|_{W_n} = \overline{w}\,|_{W_n}\) since
every factor of \(\overline{w}\) is in \(R_{n+1}\) and by Lemma
\ref{SecondIndConstr}(f), all sets in \(W_nR_{n+1}\) are disjoint
from the support of \(L_{n+1}\).  This puts \(W_0v\) in \(W_nR_{n+1}
= W_0R_1f^n\).  That \(W_0V\) contains all the \(W_nR_{n+1}\) is
clear.

We consider (4).  If \(W_0v\) intersects \(W_0\), then by Lemma
\ref{SecondIndConstr}(e) and the above argument, we must have that
the form of \(v\) is as \(w\) above: a product of elements of the
various \(H_i\).  Now let \(\overline{w}\) be the image of \(w\)
under \(\phi_0:M\rightarrow R_1\).  Our argument above has
\(v|_{W_0}=\overline{w}\,|_{W_0}\).  But \(\overline{w}\) is in
\(R_1\) and \((Z, \widetilde{Y}_1, R_1, W_0)\) is a pre-wreath
structure, so \(\overline{w}\,|_{W_0} = 1|_{W_0}\).

To work on (5), we start with a \(v\in V\).  If \(v=f^nw\) with
\(w\) a product of elements in the \(H_i\) and \(n\ne0\), then
\(W_0v\) is in \(W_nR_{n+1}\) by our analysis above.  By Lemma
\ref{SecondIndConstr}(e), \(W_0\) is not in \(W_nR_{n+1}\) so \(v\)
moves \(W_0\).  If \(n=0\), then \(v\) is in some \(R_k\).  We can
take \(k\) to be the smallest \(i\) so that \(H_i\) contains a
factor of \(w\).  Now \(v\) must move some element \(X\) of
\(W_{i-1}R_i\).  But \(X\) is of the form \(W_0v'\) for
some \(v'\in V\).

For the final claim, let \(N\) be a non-trivial normal subgroup of
\(V\).  If \(N\sub M\), then Lemma \ref{ThirdIndConstr}(c) gives the
conclusion.  If not, then an element of \(N\) is of the form
\(v=f^nw\) as above with \(n\ne0\) and with \(w\) in some \(L_i\).
Now if \(h\) is in \(H_{i-n}\), then \[\begin{split}c=v^{-1}hvh^{-1}
&= w^{-1}f^{-n}hf^nwh^{-1} \\ &= w^{-1}(h)^{f^n}wh^{-1}
\end{split}\] is in \(M\cap N\) and the calculation shows that
\(c\phi_{n-1}\) is non-trivial no matter whether \(n\) is positive
or negative, where \(\phi_{n-1}:M\rightarrow R_n\) is the
epimorphism from Lemma \ref{TopGpHom}.  Now Lemma
\ref{ThirdIndConstr}(c) applies once more.  \end{proof}

\subsection{Remarks}

There is very little in this section that is new.  As mentioned in
\cite{MR41:1884}, wreath products are often used to supply examples
of groups with special properties.  In \cite{MR41:3350}, they are
used to build groups in which others embed.  By contrast, more of
our effort goes towards recognition than creation, and so our view
of the wreath product is more ``internal'' than ``external.''
However this shift of view contains no surprises.

The reference \cite{drobinson} that we give for the wreath product
of permutation groups is not the oldest and the notion seems to have
been introduced in \cite{MR14:242b}.  Iterated wreath products
abound.  The group we call \(M\) in Lemma \ref{ThirdIndConstr} is
the group \(W\) of (5) on Page 171 of \cite{MR25:3080}.  The fact
that groups of homeomorphisms of the reals interact strongly with
wreath products is in \cite{MR41:1884} and \cite{MR41:3350} with the
latter having more overlap with the current paper.  The remarks at
the beginning of Section 5 of \cite{MR41:3350} relate to the effort
that we put into finding ``carriers'' of the various groups.

\section{Elementary amenability classes}\mylabel{ElemFactsSec}

\subsection{Squares and higher products}

We start with some elementary observations on elementary classes.
Since \(EG_0\) contains both finite groups (which might not be
abelian) and abelian groups (which might not be finite), it is clear
that \(G\times H\) might not be in \(EG_0\) even if both \(G\) and
\(H\) are in \(EG_0\).  However, \(G\times G\) will always be in
\(EG_0\) if \(G\) is in \(EG_0\).  The following can be shown as an
exercise.

\begin{lemma}\mylabel{SquaresPersist} If \(G\) is in \(EG_\alpha\),
then \(G\times G\) is also in \(EG_\alpha\).  \end{lemma}

On the other hand, if \(G\in EG_0\) is finite, then \(\Sigma(G)\),
the direct sum of countably many copies of \(G\) (all sequences in
\(G\) that are eventually the identity under pointwise
multiplication) is not in \(EG_0\).  However, if \(G\in EG_0\) is
abelian, then \(\Sigma(G)\) is in \(EG_0\).  This behavior is
important to us and we make it a definition.

\begin{defn}\mylabel{SumPersists} If \(G\in EG\) has elementary
class \(\alpha\), then we say that \(G\) has property \(\Sigma\) if
\(\Sigma(G)\) also has elementary class \(\alpha\).  \end{defn}

\subsection{Property \protect\(\Sigma\protect\) and wreath products}

There are groups of class 0 with property \(\Sigma\) because of the
abelian groups.  We get more from the next lemma.

\begin{lemma}\mylabel{WreathUp} If \((G, Z)\) is a transitive,
countably infinite, finitely generated permutation group so that
\(G\) is in \(EG\) with elementary class \(\alpha\) and with
property \(\Sigma\), then \(G\wr G\) has elementary class
\(\alpha+1\) and has property \(\Sigma\).  \end{lemma}

\begin{proof} Non-trivial wreath products are not abelian and our
hypotheses imply that \(G\wr G\) is not finite, so \(G\wr G\) does
not have class 0.  Now \(G\wr G\) is finitely generated and any
representation of \(G\wr G\) is a union of groups must have \(G\wr
G\) as one of the groups.  If \(G\wr G\) is the middle term of a short
exact sequence, then we can apply Lemma \ref{KernelControl} since
there is no problem realizing \(G\wr G\) as a result of Proposition
\ref{PreWreathToWreath}.  Now Lemma \ref{KernelControl} and
Proposition \ref{chou} say that either the kernel or quotient in the
short exact sequence has class at least that of \(G\).  Thus the
class of \(G\wr G\) is at least \(\alpha+1\).  However, by the
structure of a wreath product and the countability of \(G\), there
is a normal subgroup \(N\) of \(G\) isomorphic to \(\Sigma(G)\)
whose corresponding quotient is \(G\).  By hypothesis, the class of
\(\Sigma(G)\) is \(\alpha\), so the class of \(G\wr G\) is no more
than \(\alpha+1\).

It is not hard to represent a countable direct sum \(H\) of copies
of \(G\wr G\) as a wreath product.  The top group of \(H\) will be a
copy of \(\Sigma(G)\).  The base group of \(H\) is a countable
direct sum of copies of the base group \(N\) in \(G\wr G\).  Thus
the base group of \(H\) is isomorphic to \(\Sigma(G)\), and we get
that \(H\) ends up with class \(\alpha+1\) as well.  \end{proof}

We make the next lemma more specific to the constructions in Section
\ref{ConstrSec}.

\begin{lemma}\mylabel{WreathWayUp} Let \(\alpha\) be an ordinal and
let \(\beta\) be the smallest limit ordinal that is greater than
\(\alpha\).  Assume the notation and hypotheses of Proposition
\ref{FourthIndConstr} and assume that \((H_1, Y_0H_1)\) is a
transitive, countably infinite, finitely generated permutation group
so that \(H_1\) is in \(EG\) with elementary class \(\alpha\) and
with property \(\Sigma\).  Then the elementary class of \(M\) is
\(\beta+1\), the elementary class of \(V\) is \(\beta+2\) and both
\(M\) and \(V\) have property \(\Sigma\).  \end{lemma}

\begin{proof} All the \(H_i\) are isomorphic and the normal subgroup
\(M\) of Proposition \ref{FourthIndConstr} contains groups that can
be written as arbitrarily long sequences of the form \(H_1\wr H_1\wr
\cdots \wr H_1\).  Thus \(M\) contains all subgroups in the
inductively defined sequence \(P_0=H_1\), \(P_{i+1}=P_i\wr P_i\).
By Lemma \ref{WreathUp}, \(M\) contains subgroups of elementary
class \(\alpha+i\) for all \(i\in\N\).  Thus \(M\) is in none of the
\(EG_{\alpha+i}\) and is thus not in \(EG_\beta\) which is the union
of all the \(EG_{\alpha+i}\) by definition.  However, \(M\) is the
union of the \(P_i\) and so is in \(EG_{\beta+1}\) and the
elementary class of \(M\) is \(\beta+1\).

Since \(V=\langle H_1, f\rangle\), it is finitely generated and
getting \(V\) as a union of groups of lower class is ruled out.
Since any non-trivial normal subgroup of \(V\) contains a copy of
\(L_0\) which in turn contains a copy of every \(P_i\), the class of
\(V\) is at least \(\beta+2\).  However, \(V\) is the semidirect
product of \(M\) and a copy of the integers, so the class of \(V\)
is no more than \(\beta+2\).

Each \(P_i\) has property \(\Sigma\) and \(\Sigma(M)\) is the union
of the \(\Sigma(P_i)\), so \(M\) has property \(\Sigma\).  Now
\(\Sigma(\Z)\) has class 0, and \(\Sigma(V)\) is the semidirect
product of \(\Sigma(M)\) with \(\Sigma(\Z)\).  If follows that the
class of \(\Sigma(V)\) is that of \(V\).  \end{proof}

The next is an easy exercise from the definitions of the
\(EG_\alpha\).

\begin{lemma}\mylabel{UnionClass} If \(G_i\) is a sequence of groups
in \(EG\) whose sequence of elementary classes has no largest
element, and \(beta\) is the limit of the classes, then the direct
sum of the \(G_i\) has elementary class \(\beta+1\).  \end{lemma}

\section{The main theorem}\mylabel{MainThmSec}

\subsection{\protect\(F\protect\) or \protect\(PL_o(I)\protect\)} 

In this section we will build subgroups of \(F\), but we will spend
little time talking about \(F\).  Our groups will be groups of
homeomorphisms from the unit interval to itself and we will build
the groups by picking out homeomorphisms.  Technically, \(F\) is the
subgroup of \(PL_o(I)\) in which all slopes used are integral powers
of two, and all discontinuities of slope occur at points in
\(\Z[\frac{1}{2}]\), the set of rational numbers whose
denominators are integral powers of two \cite[Section 1]{CFP}.
Thus, we should appear to be careful and only use homeomorphisms
possessed of the properties of homeomorphisms in \(F\).  However, we
will not do so.

There is incredible flexibility possessed by elements of \(F\) that
will allow us to assume (without explicitly saying so beyond this
paragraph) that all that we build can be built in \(F\).  The key
fact is the following.  Given any two finite sequences
\(x_1<x_2<\cdots<x_k\) and \(y_1<y_2<\cdots<y_k\) of the same length
in \(\Z[\frac{1}{2}]\,\,\cap\,\,[0,1]\), there is an element of
\(F\) taking each \(x_i\) to \(y_i\).  In particular, if \(x_1=y_1\)
and \(x_k=y_k\), then the element of \(F\) can be chosen to be the
identity outside of the interval \([x_1,x_k]\).  This is recorded in
several places, such as Sections 1 and 2.2 of \cite{brin+fer}.  From
this point on, we will assume that any element called for with
certain properties in \(PL_o(I)\) can be chosen to come from \(F\).

\subsection{Bumps}

We adapt terminology from ordered groups.  All functions discussed
here are elements of \(H_o(I)\), the group of orientation preserving
self homeomorphisms of the unit interval \(I=[0,1]\).

If \(h\in H_o(I)\), then a {\itshape bump interval} of \(h\) is an
interval \([a,b]\) in \(I\) with \(a<b\) so that \(a\) and \(b\) are
fixed points of \(h\) and no \(x\in(a,b)\) is fixed by \(h\).  That
is, a bump interval of \(h\) is a maximal closed interval whose
non-empty interior contains no fixed points of \(h\).  A {\itshape
one bump function} is a function that has only one bump interval.
Two examples of one bump functions are pictured below.  

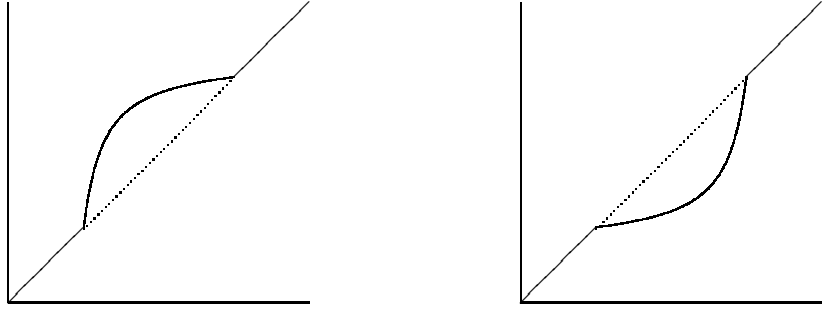
\[
\xy
(0,0); (0,40)**@{-}; (0,0); (40,0)**@{-};
(0,0); (10,10)**@{-}; (30,30)**\crv{(12,25)&(15,28)}; 
(10,10); (30,30)**@{.};
(40,40)**@{-};
\endxy
\qquad\qquad
\qquad\qquad
\xy
(0,0); (0,40)**@{-}; (0,0); (40,0)**@{-};
(0,0); (10,10)**@{-}; (30,30)**\crv{(25,12)&(28,15)}; 
(10,10); (30,30)**@{.};
(40,40)**@{-};
\endxy
\]

The word ``bump'' by itself will be dealt with less precisely.  If
\(h\) is in \(H_o(I)\), then a {\itshape bump} of \(h\) is either
the restriction of \(h\) to one of its bump intervals, or it is the
unique one bump function that agrees with \(h\) on a given bump
interval.  However, we will sometimes be careless and use the word
bump to refer to a bump interval, if at the same time it is also
being used to refer to a function on that interval.  This allows us
to make such statements as ``disjoint bumps commute.''  We hope the
reader is comfortable with this.

Given a one bump function \(h\) with bump interval \([a,b]\) a
{\itshape fundamental domain} for \(h\) is an interval \([c,d]\) in
\((a,b)\) with \(c<d\) so that \(ch=d\) or \(dh=c\).  Examples are
pictured below.

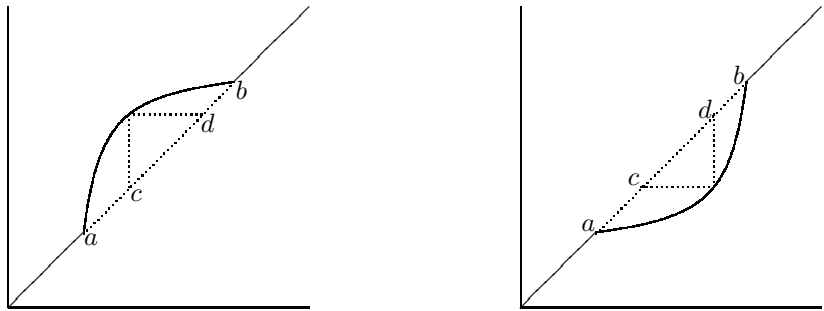
\[
\xy
(0,0); (0,40)**@{-}; (0,0); (40,0)**@{-};
(0,0); (10,10)**@{-}; (30,30)**\crv{(12,25)&(15,28)}; 
(10,10); (30,30)**@{.};
(40,40)**@{-};
(11,9)*{a}; (31,29)*{b};
(16,16); (16,25.7)**@{.}; (25.7,25.7)**@{.};
(17,15)*{c}; (26.5, 24.5)*{d};
\endxy
\qquad\qquad
\qquad\qquad
\xy
(0,0); (0,40)**@{-}; (0,0); (40,0)**@{-};
(0,0); (10,10)**@{-}; (30,30)**\crv{(25,12)&(28,15)}; 
(10,10); (30,30)**@{.};
(40,40)**@{-};
(16,16); (25.7,16)**@{.}; (25.7,25.7)**@{.};
(9,11)*{a}; (29,31)*{b};
(15,17)*{c}; (24.5, 26.5)*{d};
\endxy
\]

The extra generality in the next lemma is deliberate.

\begin{lemma}\mylabel{OneBump} If \(I_i\), \(i=1,2,3,4\), are
intervals in \(I=[0,1]\) with \(I_i\) contained in the interior of
\(I_{i+1}\), so that \(h\) is a one bump function with bump interval
\(I_3\) and fundamental domain \(I_2\) then \((I, I_4, \langle
h\rangle, I_1)\) is a pre-wreath structure.  \end{lemma}

\subsection{Finding room}

There is a lot of room in the unit interval.  We will not say a lot
about where to put the elements and groups that we create.  Room can
always be found and we will leave it to the reader to find the
room.  A key tool in finding the room is the next standard lemma.

\begin{lemma} Any countable ordinal \(\alpha\) can be embedded in an
order preserving way in \(I=[0,1]\).  \end{lemma}

\begin{proof} Let \(f\) be a bijection from \(\alpha\) to the set
\[\left\{I_i=\left[\frac{1}{2^{i+1}} , \frac{1}{2^i}\right] \mid
i\in\N\right\}.\] The embedding sends \(\beta\in\alpha\) to the
measure of the set
\(\displaystyle{\bigcup_{\gamma<\beta}{f(\gamma)},}\) and
elementary measure theory establishes the needed properties of the
embedding.  \end{proof}

\subsection{The constructions}

We start with the group generated by a one bump function.  This is
isomorphic to \(\Z\), has elementary class 0 and has property
\(\Sigma\).  This can be put in a pre-wreath setting by Lemma
\ref{OneBump}.

The inductive step is to take a finitely generated group of class
\(\alpha\) with property \(\Sigma\) in a pre-wreath structure \((I,
J, G, K)\) where \(J\) and \(K\) are intervals in the interior of
\(I\) and where the support of \(G\) has closure in the interior of
\(J\).  Then a function \(f\) can be found that takes \(K\) to
\(J\).  We can further ask that the support of \(f\) have closure in
the interior of \(I\).  Now Lemmas \ref{SecondIndConstr},
\ref{ThirdIndConstr} and Proposition \ref{FourthIndConstr} apply as
well as Lemmas \ref{WreathUp} and \ref{WreathWayUp} to give groups
with classes \(\alpha+i\) for \(i\in\N\) and \(\beta+1\) and
\(\beta+2\) where \(\beta\) is the limit of the \(\alpha+i\).  If we
call the largest group \(V(G)\), then we also get a pre-wreath
structure for \(V(G)\).  The group \(V(G)\) will be finitely
generated as well.  This process can then be repeated.

If we start with \(G_0\) as the group generated by a one-bump
function, then we can create \(G_{i+1}\) as \(V(G_i)\).  In this way
we get  groups of all classes except limit classes from class 0 up
to but not including the class \(\omega^2\).  For the class
\(\omega^2+1\), we take the direct sum of the \(G_i\) and apply
Lemma \ref{UnionClass}.  The direct sum is obtained by finding a
countable set of pairwise disjoint closed intervals in \(I\) and
conjugating each \(G_i\) into its own private interval.  The group
generated by the images of the \(G_i\) is the desired group.

We have shown Theorem \ref{MainThm}.

\subsection{Going farther}

If we let \(\widetilde{G}\) be the direct sum of the \(G_i\) of the
previous section, then \(\widetilde{G}\) is not finitely generated.
Various constructions on \(\widetilde{G}\) do not seem to raise the
class.  It is not hard to see that \(\widetilde{G}\wr
\widetilde{G}\) and \(V(\widetilde{G})\) have the same class as
\(\widetilde{G}\).

There remains the possibility that there are other construction
techniques to try.  A few attempts at other constructions gave
groups not in the class \(EG\).  In particular, they gave groups
that had subgroups isomorphic to the non-elementary amenable group
\(F\).  The presence of the group \(F\) is not hard to detect using
the main result of \cite{brin:ubiq}.

\section{Two non-solvable groups}\mylabel{ZInfinitySec}

In this section we consider the two groups \[\begin{split}
\RZ &= \Z\wr\Z\wr\Z\wr\cdots, \\ \LZ &= \cdots
\wr\Z\wr\Z\wr\Z, \end{split}\] where parentheses have been left out
because the wreath products are all wreath products of permutation
groups. 

These groups are subgroups of the group \(G_1\) of the previous
section.  The group \(G_1\) can be viewed as generated by the
functions \(h\) and \(f\) in the diagram below.  We show smooth
bumps, but the reader should interpret the pictures as representing
PL functions (or more restrictively, elements of \(F\)).

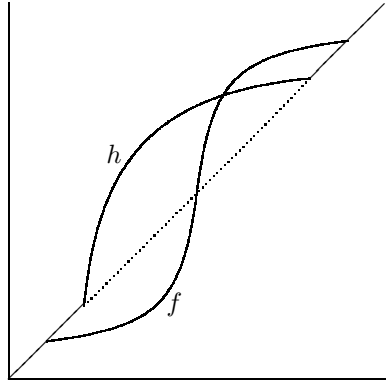
\[\xy
(0,0); (0,50)**@{-}; (0,0); (50,0)**@{-};
(0,0); (5,5)**@{-}; (45,45); (50,50)**@{-};
(5,5); (25,25)**\crv{(20,7)&(23,10)}; (45,45)**\crv{(27,40)&(30,43)};
(5,5); (10,10)**@{-}; (40,40)**\crv{(12,30)&(20,38)}; (45,45)**@{-};
(10,10); (40,40)**@{.};
(14,30)*{h}; (22,10)*{f};
\endxy
\]

The next picture shows how \(f^{-1}\) conjugates \(h\) into a one
bump function whose one bump interval fits properly inside a
fundamental domain of \(h\).

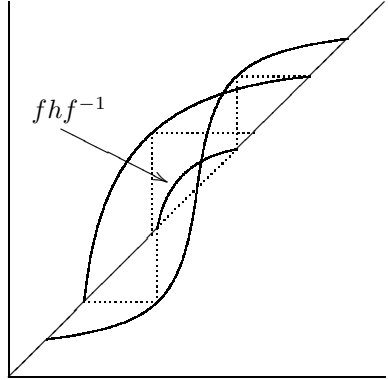
\[\xy
(0,0); (0,50)**@{-}; (0,0); (50,0)**@{-};
(0,0); (5,5)**@{-}; (45,45); (50,50)**@{-};
(5,5); (25,25)**\crv{(20,7)&(23,10)}; (45,45)**\crv{(27,40)&(30,43)};
(5,5); (10,10)**@{-}; (40,40)**\crv{(12,30)&(20,38)}; (45,45)**@{-};
(10,10); (19.7,19.7)**@{-}; (30.3,30.3)**@{.}; (40,40)**@{-};
(10,10); (19.7,10)**@{.}; (19.7, 19.7)**@{.};
(30.3,30.3)**\crv{(20.7,25.7)&(24.3,29.3)};
(30.3,40)**@{.}; (40,40)**@{.};
(7,33); (21,26)**@{-}; (21,26)*\dir{>};
(8,35.5)*{fhf^{-1}};
(19,19); (19,32.5)**@{.}; (32.5,32.5)**@{.};
\endxy
\]

If we let \(h_i=f^{-i}hf^i\) for \(i\in\Z\), then \(h=h_0\) and the
function \(fhf^{-1}\) shown is \(h_{-1}\).  It is easier to show
than \(h_1\).

The group \(\RZ\) is generated by \(h_1, h_2, \ldots\)
while \(h_{-1}, h_{-2}, \ldots\) generates \(\LZ\).

In view of the fact that we already know that \(\RZ\) and \(\LZ\)
embed in \(F\) and \(PL_o(I)\), Theorem \ref{NoCrossEmbds} will
follow from the next two lemmas.  We start with the easier of the
two.

\begin{lemma}\mylabel{NoRightInLeft} There is no subgroup of
\(\RZ\) that is isomorphic to \(\LZ\).
\end{lemma}

\begin{proof} The subgroups \(A_i=\langle h_1, h_2, \ldots,
h_i\rangle\) exhaust \(\RZ\).  If we define the groups
\(B_i=\langle h_{i+1}, h_{i+2}, \ldots\rangle\), then we have
\(\RZ=A_i\wr B_i\).  The kernel \(K_i\) of the surjection
\(\RZ\rightarrow B_i\) of Lemma \ref{TopGpHom} is
isomorphic to \(\Sigma(A_i)\), the countable direct sum of copies of
\(A_i\).  But \(A_i\) has class no greater than \(i\) (it is closer
to \(\log_2(i)\)) and \(A_i\) has property \(\Sigma\) (exercise).
Thus the class of each \(K_i\) is finite, and the \(K_i\) exhaust
\(\RZ\).

If the lemma is false, then the intersections of the \(K_i\) with
the embedded \(\LZ\) will exhaust \(\LZ\) with
normal subgroups.  But from Lemma \ref{ThirdIndConstr}, every
non-trivial normal subgroup of \(\LZ\) contains a
subgroup that surjects onto \(\LZ\) and so has class at
least that of \(\LZ\).  But \(\LZ\) contains
arbitrarily long \(\Z\wr\Z\wr\cdots\wr\Z\) iterations of wreath
products of \(\Z\) and so does not have finite class.  However, this
would have to be a subgroup of some \(K_i\) with its finite
class---a contradiction.  \end{proof}

\begin{lemma}\mylabel{NoLeftInRight} There is no subgroup of
\(\LZ\) that is isomorphic to \(\RZ\).
\end{lemma}

The following proof can be reworded in terms of the action of
\(\LZ\) on a rooted tree of infinitely many levels in which each
node has a \(\Z\) indexed set of children.  We did not feel that it
was worth making the translation.

\begin{proof} Assume a monomorphism from \(\RZ\) into
\(\LZ\).  For an element \(x\in \RZ\), we will
write its image as \(\overline{x}\).

The group \(\RZ\) is generated by \(h_1, h_2, \ldots\)
which satisfy the following whenever \(1\le i<j<k\):
\mymargin{TheLaw}\begin{equation}\label{TheLaw}\tag{\protect\(*\protect\)}
\left[ h_i\,\,,\,\, (h_{k})^{-n}(h_{j})(h_{k})^n \right] = 1
\quad\hbox{if and only if}\quad n\ne0.
\end{equation}
In words, \(h_{j}\) is the only conjugate of \(h_{j}\) by a
power of \(h_{k}\) that fails to commute with \(h_i\).

We will focus on the fact that \tref{TheLaw} must be satisfied by
the \(\overline{h}_1, \overline{h}_2, \ldots\) in
\(\LZ\).  

We will talk a great deal about bumps.

Elements of \(\LZ\) are made of finitely many pairwise disjoint
bumps.  Given two bumps from different elements of \(\LZ\), the
intervals are either disjoint, identical, or one is in the interior
of the other.  This is not hard to show as an exercise (which can
further show that bump intervals in \(\LZ\) are very restricted in
that if the left endpoint of a bump is known, then its right
endpoint is determined and vice versa).

If two bumps are related by having the interval of one in the
interior of the interval of the other, then we will say that they
are ``nested.''  The bump with the larger interval will be called
``superior to'' the bump with the smaller interval and the inverse
relation will be called ``inferior to.''

It is extremely important that an inferior bump in \(\LZ\) is
always contained in a fundamental domain of a superior bump.

Nested bumps cannot commute.

Another important fact that we need is that in \(\LZ\), any bump
interval is inferior to only finitely many other bump intervals from
\(\LZ\).

We need a sublemma: {\itshape If \(u\), \(v\) and \(w\) are bumps of
piecewise linear functions on the same bump interval \([a,b]\) and
if \(w^{-m}vw^m\) and \(w^{-n}vw\) both commute with \(u\), then
\(w^{n-m}\) commutes with \(v\).}

\begin{proof}[Proof of the sublemma] This is the only argument in
the paper that uses the PL setting.  A function with bump interval
\([a,b]\) that commutes with \(u\) is completely determined by what
it does on a fundamental domain of \(u\) in \([a,b]\).  There are
infinitely many fundamental domains of \(u\) in any
\([a,a+\epsilon)\) with \(\epsilon>0\).  Because the functions are
piecewise linear and \(a\) is a fixed point of \(v\) and \(w\),
there is an \([a,a+\epsilon)\) on which all of \(v\), \(w^{-m}vw^m\)
and \(w^{-n}vw^n\) agree.  Since the last two commute with \(u\),
they are identical and the sublemma is proved.  \end{proof}

\noindent{\itshape Proof of Lemma \ref{NoLeftInRight}, continued.}
We first argue that if there are \(i<j\) so that \(\overline{h}_i\)
and \(\overline{h}_j\) have nested bumps, then we have a
contradiction.  From \tref{TheLaw}, we know that
\(\overline{h}_{j+1}\) must conjugate \(\overline{h}_j\) so that no
bumps from \(\overline{h}_i\) and the conjugate of
\(\overline{h}_j\) are nested.  This can only happen if a bump of
\(\overline{h}_{j+1}\) is superior to the relevant bumps of both
\(\overline{h}_i\) and \(\overline{h}_j\).  Now the argument can be
repeated for \(\overline{h}_j\) and \(\overline{h}_{j+1}\) since
they have nested bumps.  In this way we create an infinite ascending
chain of nested bumps.  As observed above, this is not possible in
\(\LZ\).

We now assume that for all \(i\ne j\), if bump intervals of
\(\overline{h}_i\) and \(\overline{h}_j\) intersect, then they are
identical.  Since \(\overline{h}_1\) and \(\overline{h}_2\) do not
commute, there has to be at least one interval \([a,b]\) that is a
bump interval for bumps of \(\overline{h}_1\) and \(\overline{h}_2\)
that do not commute.  By \tref{TheLaw}, all non-zero powers of
\(\overline{h}_3\) conjugate \(\overline{h}_2\) so that the results
commute with \(\overline{h}_1\).  By our latest assumption, this is
accomplished by a bump of \(\overline{h}_3\) whose interval is also
\([a,b]\).  By the sublemma, \(\overline{h}_3\) commutes with
\(\overline{h}_2\) on \([a,b]\).  But this makes it impossible for
\(\overline{h}_3\) to alter \(\overline{h}_2\) on \([a,b]\) by
conjugation.  Since \(\overline{h}_2\) does not commute with
\(\overline{h}_1\) on \([a,b]\) we contradict \tref{TheLaw}.  This
proves the lemma.  \end{proof}


\begin{thebibliography}{10}

\bibitem{MR84m:43001}
S.~I. Adyan, \emph{Random walks on free periodic groups}, Izv. Akad. Nauk SSSR
  Ser. Mat. \textbf{46} (1982), no.~6, 1139--1149, 1343. \MR{84m:43001}

\bibitem{brin:ubiq}
Matthew~G. Brin, \emph{The ubiquity of {T}hompson's group ${F}$ in groups of
  piecewise linear homeomorphisms of the unit interval}, J. London Math. Soc.
  (2) \textbf{60} (1999), no.~2, 449--460. \MR{1 724 861}

\bibitem{brin+fer}
Matthew~G. Brin and Fernando Guzm{\'a}n, \emph{Automorphisms of generalized
  {T}hompson groups}, J. Algebra \textbf{203} (1998), no.~1, 285--348.

\bibitem{brin+squier}
Matthew~G. Brin and Craig~C. Squier, \emph{Groups of piecewise linear
  homeomorphisms of the real line}, Invent. Math. \textbf{79} (1985), 485--498.

\bibitem{brown:finiteprop}
Kenneth~S. Brown, \emph{Finiteness properties of groups}, Journal of Pure and
  Applied Algebra \textbf{44} (1987), 45--75.

\bibitem{MR2001k:20087}
J.~Burillo, S.~Cleary, and M.~I. Stein, \emph{Metrics and embeddings of
  generalizations of {T}hompson's group {$F$}}, Trans. Amer. Math. Soc.
  \textbf{353} (2001), no.~4, 1677--1689 (electronic). \MR{2001k:20087}

\bibitem{burillo}
Jos{\'e} Burillo, \emph{Quasi-isometrically embedded subgroups of {T}hompson's
  group ${F}$}, J. Algebra \textbf{212} (1999), no.~1, 65--78.

\bibitem{CFP}
J.~W. Cannon, W.~J. Floyd, and W.~R. Parry, \emph{Introductory notes on
  {R}ichard {T}hompson's groups}, Enseign. Math. (2) \textbf{42} (1996),
  no.~3-4, 215--256. \MR{98g:20058}

\bibitem{chou}
Ching Chou, \emph{Elementary amenable groups}, Illinois J. Math. \textbf{24}
  (1980), no.~3, 396--407. \MR{81h:43004}

\bibitem{MR2016187}
Sean Cleary and Jennifer Taback, \emph{Geometric quasi-isometric embeddings
  into {T}hompson's group {$F$}}, New York J. Math. \textbf{9} (2003), 141--148
  (electronic). \MR{2 016 187}

\bibitem{MR19:1067c}
Mahlon~M. Day, \emph{Amenable semigroups}, Illinois J. Math. \textbf{1} (1957),
  509--544. \MR{19,1067c}

\bibitem{dehornoy:assoc}
Patrick Dehornoy, \emph{The structure group for the associativity identity}, J.
  Pure Appl. Algebra \textbf{111} (1996), no.~1-3, 59--82. \MR{97d:55033}

\bibitem{dydak2}
J.~Dydak, \emph{1-movable continua need not be pointed 1-movable}, Bull. Acad.
  Polon. Sci. S{\'e}r. Sci. Math. Astronom. Phys. \textbf{25} (1977), 485--488.

\bibitem{dydak}
\bysame, \emph{A simple proof that pointed, connected {F}{A}{N}{R}-spaces are
  regular fundamental retracts of {A}{N}{R}'s}, Bull. Acad. Polon. Sci.
  S{\'e}r. Sci. Math. Astronom. Phys. \textbf{25} (1977), 55--62.

\bibitem{freyd+heller}
Peter Freyd and Alex Heller, \emph{Splitting homotopy idempotents {I}{I}},
  Journal of Pure and Applied Algebra \textbf{89} (1993), 93--106.

\bibitem{gersten:problems}
S.~M. Gersten, \emph{Selected problems}, Combinatorial Group Theory and
  Topology (S.~M. Gersten and John~R. Stallings, eds.), Annals of Mathematics
  Studies, vol. 111, Princeton {U}niversity {P}ress, 1987, pp.~545--551.

\bibitem{grig2}
R.~I. Grigorchuk, \emph{Degrees of growth of finitely generated groups and the
  theory of invariant means}, Izv. Akad. Nauk SSSR Ser. Mat. \textbf{48}
  (1984), no.~5, 939--985. \MR{86h:20041}

\bibitem{MR99b:20055}
\bysame, \emph{An example of a finitely presented amenable group that does not
  belong to the class {EG}}, Mat. Sb. \textbf{189} (1998), no.~1, 79--100.
  \MR{99b:20055}

\bibitem{guba+sapir:diag2}
V.~S. Guba and M.~V. Sapir, \emph{On subgroups of the {R}. {T}hompson group
  {$F$} and other diagram groups}, Mat. Sb. \textbf{190} (1999), no.~8, 3--60.
  \MR{2001m:20045}

\bibitem{guba+sapir2}
Victor Guba and Mark Sapir, \emph{Diagram groups}, Mem. Amer. Math. Soc.
  \textbf{130} (1997), no.~620, viii+117. \MR{98f:20013}

\bibitem{MR25:3080}
P.~Hall, \emph{Wreath powers and characteristically simple groups}, Proc.
  Cambridge Philos. Soc. \textbf{58} (1962), 170--184. \MR{25 \#3080}

\bibitem{MR41:1884}
W.~Charles Holland, \emph{The characterization of generalized wreath products},
  J. Algebra \textbf{13} (1969), 152--172. \MR{41 \#1884}

\bibitem{MR41:3350}
W.~Charles Holland and Stephen~H. McCleary, \emph{Wreath products of ordered
  permutation groups}, Pacific J. Math. \textbf{31} (1969), 703--716. \MR{41
  \#3350}

\bibitem{MR14:242b}
L{\'e}o Kaloujnine and Marc Krasner, \emph{Produit complet des groupes de
  permutations et probl\`eme d'extension de groupes. {I}}, Acta Sci. Math.
  Szeged \textbf{13} (1950), 208--230. \MR{14,242b}

\bibitem{neumann:wreath}
Peter~M. Neumann, \emph{On the structure of standard wreath products of
  groups}, Math. Z. \textbf{84} (1964), 343--373. \MR{32 \#5719}

\bibitem{MR82b:43002}
A.~Ju. Ol$'${\v{s}}anski\u{\i}, \emph{On the question of the existence of an
  invariant mean on a group}, Uspekhi Mat. Nauk \textbf{35} (1980), no.~4(214),
  199--200. \MR{82b:43002}

\bibitem{MR2004f:20061}
Alexander~Yu. Ol$'$shanskii and Mark~V. Sapir, \emph{Non-amenable finitely
  presented torsion-by-cyclic groups}, Publ. Math. Inst. Hautes \'Etudes Sci.
  (2002), no.~96, 43--169 (2003). \MR{2004f:20061}

\bibitem{drobinson}
Derek J.~S. Robinson, \emph{A course in the theory of groups}, 2nd ed.,
  Graduate Texts in Math., vol.~80, Springer, New York, 1996.

\bibitem{thomp:notes}
Richard~J. Thompson, 1973, Handwritten, widely circulated, unpublished notes
  attributed to Thompson, (c. 1973+).

\bibitem{MR87e:04007}
Stan Wagon, \emph{The {B}anach-{T}arski paradox}, Encyclopedia of Mathematics
  and its Applications, vol.~24, Cambridge University Press, Cambridge, 1985,
  With a foreword by Jan Mycielski. \MR{87e:04007}

\bibitem{MR94g:04005}
\bysame, \emph{The {B}anach-{T}arski paradox}, Cambridge University Press,
  Cambridge, 1993, With a foreword by Jan Mycielski, Corrected reprint of the
  1985 original. \MR{94g:04005}

\end{thebibliography}

\providecommand{\bysame}{\leavevmode\hbox to3em{\hrulefill}\thinspace}
\providecommand{\MR}{\relax\ifhmode\unskip\space\fi MR }
\providecommand{\MRhref}[2]{%
  \href{http://www.ams.org/mathscinet-getitem?mr=#1}{#2}
}
\providecommand{\href}[2]{#2}

\noindent Department of Mathematical Sciences

\noindent State University of New York at Binghamton

\noindent Binghamton, NY 13902-6000

\noindent USA

\noindent email: matt@math.binghamton.edu

\end{document}